\documentclass[leqno,12pt]{article}

\usepackage{amssymb}
\usepackage{euscript}
\usepackage{mathrsfs} 
\usepackage[dvips]{graphicx}
\usepackage{flafter}
\usepackage[latin1]{inputenc}
\usepackage{amsmath}
\usepackage{amsthm}
\usepackage{amsfonts}
\usepackage{hyperref} 
\usepackage{amsxtra} 
\usepackage{verbatim} 
\usepackage{color} 

\evensidemargin\oddsidemargin

\def\r{{\mathbb R}}

\def\P{{\bf P}}
\def\E{{\bf E}}

\def\z{{\mathbb Z}}

\def\d{\, \mathrm{d}}

\def\DD{{\mathcal D}}

\begin{document}


\newtheorem{theorem}{Theorem}
\newtheorem{lemma}{Lemma}
\newtheorem{proposition}{Proposition}
\newtheorem{corollary}{Corollary}
\newtheorem{question}{Question}
\newtheorem{conjecture}{Conjecture}

\theoremstyle{definition}   
\newtheorem{remark}{Remark}
\newtheorem{example}{Example}

\baselineskip=18pt
\setcounter{page}{1}


\vglue50pt

\centerline{\large\bf A max-type recursive model: some properties}

\medskip

\centerline{\large\bf and open questions}

\bigskip
\bigskip

\centerline{Xinxing Chen, Bernard Derrida, Yueyun Hu, Mikhail Lifshits and Zhan Shi}

\bigskip
\bigskip

\centerline{\it In honor of Professor Charles M.\ Newman}

\centerline{\it on the occasion of his 70th birthday}

\bigskip
\bigskip
\bigskip

{\leftskip=1.5truecm \rightskip=1.5truecm \baselineskip=15pt \small

\noindent{\slshape\bfseries Summary.} We consider a simple max-type recursive model which was introduced in the study of depinning transition in presence of strong disorder, by Derrida and Retaux~\cite{derrida-retaux}. Our interest is focused on the critical regime, for which we study the extinction probability, the first moment and the moment generating function. Several stronger assertions are stated as conjectures.

\bigskip


\bigskip


} 

\bigskip
\bigskip

\section{Introduction}
\label{s:intro}

Fix an integer $m\ge 2$. Let $X_0 \ge 0$ be a non-negative random variable (but very soon, taking values in $\z_+ := \{ 0, \, 1, \, 2, \ldots\}$). To avoid trivial discussions, we assume that $X_0$ is not almost surely a constant. Consider the following recurrence relation:
\begin{equation}
    X_{n+1}
    =
    (X_n^{(1)} + \cdots + X_n^{(m)} -1)^+ ,
    \qquad
    n\ge 0 ,
    \label{iter}
\end{equation}

\noindent where $X_n^{(1)}$, $\ldots$, $X_n^{(m)}$ are independent copies of $X_n$, and for all $x\in \r$, $x^+ := \max\{ x, \, 0\}$ is the positive part of $x$.

The model with recursion defined in \eqref{iter} was introduced by Derrida and Retaux~\cite{derrida-retaux} as a simple hierarchical renormalization model to understand depinning transition of a line in presence of strong disorder. The study of depinning transition has an important literature both in mathematics and in physics. Of particular interest are problems about the relevance of the disorder, and if it is, the precise description of the transition. Similar problems are raised when the line has hierarchical constraints. We refer to Derrida, Hakim and Vannimenus~\cite{derrida-hakim-vannimenus}, Giacomin, Lacoin and Toninelli~\cite{giacomin-lacoin-toninelli}, Lacoin~\cite{lacoin}, Berger and Toninelli~\cite{berger-toninelli}, Derrida and Retaux~\cite{derrida-retaux}, and Hu and Shi~\cite{yz_bnyz} for more details and references. Let us mention that the recursion \eqref{iter} appears as a special case in the survey paper of Aldous and Bandyopadhyay~\cite{aldous-bandyopadhyay}, in a spin glass toy-model in Collet et al.\ \cite{collet-eckmann-glaser-martin}, and is also connected to a parking scheme recently studied by Goldschmidt and Przykucki~\cite{goldschmidt-przykucki}.

Since $x-1 \le (x-1)^+ \le x$ for all $x\ge 0$, we have, by \eqref{iter},
$$
   m \, \E(X_n) -1 \le \E(X_{n+1}) \le m \, \E(X_n),
$$

\noindent so the {\it free energy}
\begin{equation}
    F_\infty
    :=
    \lim_{n\to \infty} \downarrow \, \frac{\E(X_n)}{m^n}
    =
    \lim_{n\to \infty} \uparrow \, \frac{\E(X_n) - \frac{1}{m-1}}{m^n}\, ,
    \label{F}
\end{equation}

\noindent is well-defined.

In a sense, the free energy is positive for ``large'' random variables and vanishes for ``small'' ones. The two regimes are separated by a ``surface'' (called ``critical manifold" in \cite{derrida-retaux}) in the space of distributions that exhibits a critical behavior. Many interesting questions are related to the behavior of $X_n$ at the critical regime or near it.

As an example, let us recall a key conjecture, due to Derrida and Retaux~\cite{derrida-retaux}, which handles the following parametric setting.

For any random variable $X$, we write $P_X$ for its law. Assume
$$
     P_{X_0}
     =
     (1-p) \, \delta_0 + p\, P_{Y_0},
$$

\noindent where $\delta_0$ is the Dirac measure at $0$, $Y_0$ is a positive random variable, and $p\in [0, \, 1]$ a parameter. Since $F_\infty =: F_\infty(p)$ is non-decreasing in $p$, there exists $p_c \in [0, \, 1]$ such that $F_\infty >0$ for $p>p_c$ and that $F_\infty (p) =0$ for $p<p_c$. A conjecture of Derrida and Retaux~\cite{derrida-retaux} says that if $p_c>0$ (and possibly under some additional integrability conditions on $Y_0$), then
\begin{equation}
     F_\infty(p)
     =
     \exp \Big( - \frac{K+o(1)}{(p-p_c)^{1/2}} \Big),
     \qquad
     p\downarrow p_c\, ,
     \label{conj:bernard}
\end{equation}

\noindent for some constant $K\in (0, \, \infty)$.

When $p_c=0$, it is possible to have other exponents than $\frac12$ in \eqref{conj:bernard}, see \cite{yz_bnyz}, which also contains several open problems in the regime $p\downarrow p_c\,$.

We have not been able to prove or disprove the conjecture. In this paper, we are interested in the critical regime, i.e., $p=p_c$ in the Derrida--Retaux conjecture setting. However, we do not formulate the model in a parametric way. When $X_0$ is integer valued, the critical regime is characterized by the following theorem.

\medskip

\begin{theorem}
\label{t:collet-et-al}

 {\bf (Collet et al.~\cite{collet-eckmann-glaser-martin}).}
 If $X_0$ takes values in $\z_+ := \{ 0, \, 1, \, 2, \ldots\}$, the critical regime is given by
 $$
     (m-1) \E(X_0 \, m^{X_0}) = \E(m^{X_0})<\infty \, ;
 $$

\noindent more precisely, $F_\infty>0$ if either $\E(m^{X_0}) < (m-1) \E(X_0 \, m^{X_0}) <\infty$, or $\E(X_0 \, m^{X_0}) =\infty$ (a fortiori, if $\E(m^{X_0}) =\infty$), and $F_\infty =0$ otherwise.

\end{theorem}

\medskip

We assume {\bf from now on} that $X_0$ is $\z_+$-valued, and we work in the {\bf critical regime}, i.e., assuming
\begin{equation}
    (m-1) \E(X_0 \, m^{X_0}) = \E(m^{X_0}) < \infty \, .
    \label{variete-critique}
\end{equation}

\noindent A natural question is whether $\E (X_n) \to 0$ in the critical regime. The answer is positive.

\medskip

\begin{theorem}
\label{t:G(Xn)->1}

 Assume $\eqref{variete-critique}$. Then $\lim_{n\to \infty} \E(m^{X_n}) =1$. A fortiori, $\lim_{n\to \infty} \E(X_n) = 0$.

\end{theorem}

\medskip

It is natural to study the asymptotic behavior of $X_n$ quantitatively.
Although  we have not succeeded in making many of our arguments rigorous,
we are led by a general asymptotic picture described by the following two
conjectures. The first of them (Conjecture \ref{conj:P(Xn=0)}) describes how $\P(X_n \not= 0)$ tends to $0$, while the second one (Conjecture \ref{conj:cvg_faible}) describes the conditional asymptotic behavior of $X_n$ provided that $X_n \not= 0$.

We use the notation $a_n \sim b_n$, $n\to \infty$, to denote $\lim_{n\to \infty} \frac{a_n}{b_n} =1$.

\medskip

\begin{conjecture}
\label{conj:P(Xn=0)}

 Assume $\eqref{variete-critique}$. Then we have
 $$
    \P (X_n \not= 0) \, \sim \, \frac{4}{(m-1)^2} \, \frac{1}{n^2},
    \qquad
    n\to \infty .
 $$

\end{conjecture}

\medskip

When $m=2$, Conjecture \ref{conj:P(Xn=0)} was already given by Collet et al.~\cite{collet-eckmann-glaser-martin} (see their equation (AIII.10)) and by Derrida and Retaux~\cite{derrida-retaux}.

Our next conjecture concerns weak convergence of $X_n$ given $X_n >0$.

\medskip

\begin{conjecture}
\label{conj:cvg_faible}

 Assume $\eqref{variete-critique}$. Then, conditionally on $X_n \not= 0$, the random variable $X_n$ converges weakly to a limit $Y_\infty$ with geometric law:
 $\P(Y_\infty = k) = \frac{m-1}{m^k}$ for integers $k\ge 1$.

\end{conjecture}

\medskip

The two conjectures above immediately lead to more specific quantitative assertions. In view of Theorem \ref{t:G(Xn)->1}, it is natural to study how $\E(X_n)$ goes to zero in the critical regime.

\begin{conjecture}
\label{conj:E(Xn)}

 Assume $\eqref{variete-critique}$. Then we have
 $$
 \E(X_n) \, \sim \, \frac{4m}{(m-1)^3} \, \frac{1}{n^2},
 \qquad
 n\to \infty \, ,
 $$

\noindent and more generally, for all real numbers $r\in (0, \, \infty)$,
 $$
 \E(X_n^r) \, \sim \, \frac{c(r)}{n^2},
 \qquad
 n\to \infty \, ,
 $$

\noindent where $c(r) = c(r, \, m) := \frac{4}{m-1} \sum_{k=1}^\infty \frac{k^r}{m^k}$.

\end{conjecture}

\medskip

In view of Conjectures \ref{conj:P(Xn=0)} and \ref{conj:cvg_faible}, we may also guess how fast the moment generating function converges.

\medskip

\begin{conjecture}
\label{conj:G(Xn)}

 Assume $\eqref{variete-critique}$. Then

{\rm a)} We have, for $n\to \infty$ and $ s\in (0, \, m), \; s\not= 1$,
 \begin{equation}
     \E(s^{X_n}) -1
     \sim
     \frac{4m}{(m-1)^2} \, \frac{s-1}{m-s} \, \frac{1}{n^2}\, .
     \label{G(X_n)(s)-1}
  \end{equation}

{\rm b)} We have, for $n\to \infty$,
  \begin{equation}
     \E(m^{X_n}) -1
     \sim
     \frac{2}{m-1} \, \frac{1}{n}\, .
     \label{G(X_n)(m)-1}
 \end{equation}

\end{conjecture}

\medskip

Possibly some additional integrability conditions, such as $\E(X_0^3 \, m^{X_0})<\infty$ in Theorem \ref{t:G(Xn)(m)-1} below, are necessary for our conjectures to hold.

The following weaker version of \eqref{G(X_n)(m)-1} can be rigorously proved.

\medskip

\begin{theorem}
\label{t:G(Xn)(m)-1}

 Assume $\eqref{variete-critique}$. If $\Lambda_0:= \E(X_0^3 \, m^{X_0})<\infty$, then there exist constants $c_2 \ge c_1>0$, depending only on $m$, such that
 $$
 c_1 \, \frac{[\P(X_0=0)]^{1/2}}{\Lambda_0^{1/2}} \, \frac1n
 \le
 \E(m^{X_n}) -1
 \le
 c_2 \, \frac{\Lambda_0^{1/2}}{[\P(X_0=0)]^{1/2}} \, \frac1n\, ,
 \qquad
 n\ge 1 \, .
 $$

\end{theorem}

\medskip

\noindent {\bf Remark.} (i) The assumption \eqref{variete-critique} guarantees that $\P(X_0=0) >0$.

(ii) It will be seen (in Lemma \ref{l:elementaire}) that in case $m\ge 3$, we have $\P(X_0=0) \ge \frac{m-2}{m-1}$.\qed

\bigskip

The content of the rest of the paper is as follows:

$\bullet$ Section \ref{s:G}: A few key facts about the moment generating functions of $X_n$ and their interrelations, with some technical proofs postponed to Section \ref{s:prop_proofs};

$\bullet$ Section \ref{s:proof_of_Theorem}: Proof of Theorem \ref{t:collet-et-al};

$\bullet$ Section \ref{s:pf_Theorem_E(Xn)->0}: Proof of Theorem \ref{t:G(Xn)->1};

$\bullet$ Section \ref{s:conj_P(Xn=0)}: Heuristics for Conjectures \ref{conj:P(Xn=0)} and \ref{conj:G(Xn)}b;

$\bullet$ Section \ref{s:conj_E(Xn)}: Heuristics for Conjectures
\ref{conj:cvg_faible},  \ref{conj:E(Xn)}, and \ref{conj:G(Xn)}a;

$\bullet$ Section \ref{s:pf_Theorem_G(Xn)(m)-1}: Proof of Theorem \ref{t:G(Xn)(m)-1};

$\bullet$ Section \ref{s:harmonic_limit}: Proofs of weaker versions of some of our conjectures.

\medskip

Throughout the paper, $c_i=c_i(m)$, for $1\le i\le 17$, denote finite and positive constants whose values depend only on $m$.

\section{Moment generating functions and their evolution}
\label{s:G}

\subsection{Derivatives and evolution}

All the techniques used in this article are based on the evaluation of the moment generating functions 
and on their evolution during the recursive process \eqref{iter}. Write, for $n\ge 0$,
$$
   G_n(s)
   :=
   \E (s^{X_n}),
$$

\noindent the moment generating function of $X_n$.

In terms of generating functions, the recursion \eqref{iter} writes as
\begin{equation}
    \label{G0}
    G_{n+1}(s)
    = \frac1s \, G_n(s)^m + (1-\frac1s)G_n(0)^m \, .
\end{equation}

\noindent Moreover, if $G_n'(s)$ is well defined, then so is $G_{n+1}'(s)$ and differentiation yields
\begin{equation}
    \label{G1}
    G_{n+1}'(s)
    = \frac{m}{s}\, G_n'(s) \, G_n(s)^{m-1} - \frac{1}{s^2} \, G_n(s)^m + \frac{1}{s^2} \, G_n(0)^m \, .
\end{equation}

\noindent Eliminating $G_n(0)$ from the two identities, it follows that
\begin{equation}
    \label{G10}
    (s-1)s\, G_{n+1}'(s) - G_{n+1}(s)
    =
    [ m(s-1) \, G_n'(s) - G_n(s) ] \, G_n(s)^{m-1} \, .
\end{equation}

\noindent Taking $s=m$ yields a formula particularly convenient for iterations,
namely,
\begin{equation}
    \label{G10(m)}
    (m-1)m\, G_{n+1}'(m) - G_{n+1}(m)
    =
    [ (m-1)m \, G_n'(m) - G_n(m) ] \, G_n(m)^{m-1} \, .
\end{equation}

Further differentiation of \eqref{G1} yields (if the involved derivatives are well defined for $G_0$)
\begin{eqnarray}
    \label{G2}
    G_{n+1}''(s)
 &=& - \frac{2m}{s^2}\, G_n'(s) \, G_n(s)^{m-1}
    +
    \frac{m}{s}\, G_n''(s) \, G_n(s)^{m-1}
    \\
    \nonumber
 &&\qquad
    +   \frac{m(m-1)}{s}\, G_n'(s)^2 \, G_n(s)^{m-2}
    +   \frac{2}{s^3} \, G_n(s)^m
    -   \frac{2}{s^3} \, G_n(0)^m
\end{eqnarray}

\noindent and
\begin{eqnarray}
    \label{G3}
 && G_{n+1}'''(s)
    = \frac{6m}{s^3}\, G_n'(s) \, G_n(s)^{m-1}
    -
    \frac{3m}{s^2}\, G_n''(s) \, G_n(s)^{m-1}
  \\   \nonumber
 &&\quad
    -
    \frac{3m(m-1)}{s^2}\, G_n'(s)^2 \, G_n(s)^{m-2}
    - \frac{6}{s^4} \, G_n(s)^m
    + \frac{6}{s^4} \, G_n(0)^m
    + \frac{m}{s} \, G_n'''(s) G_n(s)^{m-1}
    \\   \nonumber
 &&\quad
    +  \frac{3m(m-1)}{s} \, G_n'(s)G_n''(s)G_n(s)^{m-2}
    +  \frac{m(m-1)(m-2)}{s} \, G_n'(s)^3G_n(s)^{m-3} \, .
\end{eqnarray}

\noindent Notice that assumption \eqref{variete-critique} can be rewritten in the language of generating functions as
$$
    (m-1)m \, G_0'(m)= G_0(m).
$$

\noindent It follows immediately from \eqref{G10(m)} that we have in this case for all $n\ge 0$,
\begin{equation}
    \label{variete-critique-Gn}
    (m-1)m \, G_n'(m)= G_n(m).
\end{equation}

Assuming $\E(X_0^2 \, m^{X_0}) <\infty$ and plugging \eqref{variete-critique-Gn} into \eqref{G2} with $s=m$, we obtain, for all $n\ge 0$,
\begin{equation}
    \label{G2(m)}
    G_{n+1}''(m) =
        G_n''(m) \, G_n(m)^{m-1}
    +   \frac{m-2}{m^3(m-1)} \, G_n(m)^m
    -   \frac{2}{m^3} \, G_n(0)^m,
\end{equation}

\noindent which, in combination with \eqref{G0}, yields
\begin{equation}
    \label{G20}
    G_{n+1}''(m) + \frac{2}{m^2(m-1)} \, G_{n+1}(m)
    =
    \big[G_n''(m) + \frac{1} { m^2 (m-1)}\, G_n(m)\big] \, G_n(m)^{m-1} \, .
\end{equation}

\noindent This is still not a perfectly iterative relation because of the difference of the nominators on both sides.

However, let us continue the same operation with the third derivative.
Assuming that $\E(X_0^3 \, m^{X_0}) <\infty$ and plugging \eqref{variete-critique-Gn} into \eqref{G3} with $s=m$, we obtain for all $n\ge 0$,
\begin{equation}
    \label{G3(m)}
    G_{n+1}'''(m)
    =
     \frac{-2m^2+7m-6}{m^4(m-1)^2} \, G_n(m)^m
    + \frac{6}{m^4} \, G_n(0)^m
    +  G_n'''(m) G_n(m)^{m-1}.
\end{equation}

\noindent By excluding $G_n(0)$ from \eqref{G3(m)} and \eqref{G2(m)}, it follows that
\begin{eqnarray}
 &&m G_{n+1}'''(m)+3 G_{n+1}''(m)
    \nonumber
    \\
 &=& \big[ m G_n'''(m)+3 G_n''(m) \big]  G_n(m)^{m-1}
    +
    \frac{m-2}{m^2(m-1)^2} \, G_n(m)^m.
    \label{G32}
\end{eqnarray}

\noindent Now we may get a completely iterative aggregate
\begin{eqnarray}
 && \DD_n(m)
    :=
    m(m-1)G_n'''(m) + (4m-5)G_n''(m) + \frac{2(m-2)}{m^2(m-1)} \, G_n(m)
    \nonumber
    \\
 &=&
    (m-1) \big[m G_n'''(m) + 3 G_n''(m)\big]
    +
    (m-2) \big[ G_n''(m) + \frac{2}{m^2(m-1)} \, G_n(m)\big]
    \label{D2(m)}
\end{eqnarray}

\noindent because combining \eqref{G20} and \eqref{G32} yields, for all $n\ge 0$,
\begin{equation}
    \label{D2(m):lemme}
    \DD_{n+1}(m)
    =
    \DD_n(m)\, G_n(m)^{m-1}.
\end{equation}

\noindent This is another perfectly iterative relation along with \eqref{G10(m)}.
Recall that it is proved under $\eqref{variete-critique}$ and assuming that
$\E(X_0^3 \, m^{X_0}) <\infty$.

The approach based on analysis of the moment generating function is already adopted by Collet et al.~\cite{collet-eckmann-glaser-martin} and Derrida and Retaux~\cite{derrida-retaux}, where the main attention is focused on the case $m=2$. As seen from \eqref{D2(m)}, the case $m=2$ only involves the first half of $\DD_n(m)$. Our work reveals the importance of the second half; together with the first half, they serve as a useful tool in the study of the moment generating function, as we will see in Section \ref{subs:proof_Pi>}.

\subsection{Products}
\label{ss:products}

The main technical properties of generating functions are contained in the following two complementary propositions.

\medskip

\begin{proposition}
\label{p:Pi<}

 Assume $\eqref{variete-critique}$. There exists a constant $c_3 \in (0, \, \infty)$ such that
 $$
 \prod_{i=0}^n G_i(m)^{m-1}
 \le
 \frac{c_3}{\P(X_0=0)}\, n^2 \, ,
 \qquad
 n\ge 1.
 $$

\end{proposition}

\medskip

\begin{proposition}
\label{p:Pi>}

Assume $\eqref{variete-critique}$. If $\Lambda_0 := \E(X_0^3 \, m^{X_0})<\infty$, then there exists a positive constant $c_4$, depending only on $m$, such that
 $$
 \prod_{i=0}^n G_i(m)^{m-1}
    \ge
    \frac{c_4}{\Lambda_0} \, n^2 \, ,
    \qquad
    n\ge 1.
 $$

\end{proposition}

\bigskip

The proofs of these propositions, rather technical, are postponed to Section \ref{s:prop_proofs}.

\section{Proof of Theorem \ref{t:collet-et-al}}
\label{s:proof_of_Theorem}

The proof of Theorem \ref{t:collet-et-al} essentially follows Collet
et al.~\cite{collet-eckmann-glaser-martin}. It is presented here, not only for the sake of self-containedness, but also for some simplification, which we consider as interesting, in both the upper and the lower bounds.

Assume for a while that $\E(X_0 \, m^{X_0}) <\infty$, which means, in the language of generating functions, that $G_0'(s)$ is well defined for all $0\le s \le m$. Then, as we know from \eqref{G1}, the derivative $G_n'(s)$ is well defined for all $n\ge 0$ and $0\le s\le m$.
\medskip

\noindent {\it Proof of Theorem \ref{t:collet-et-al}: Upper bound.} Suppose $(m-1)\E(X_0 \, m^{X_0}) \le \E(m^{X_0})$ and $\E(X_0 \, m^{X_0})<\infty$. Let us prove that $F_\infty = 0$.
In the language of generating functions our assumption simply means $(m-1)m\, G_{0}'(m) - G_{0}(m) \le 0$ and $G_0'(m)<\infty$. Then iterative identity \eqref{G10(m)} yields that we have the same relations $(m-1)m\, G_{n}'(m) - G_{n}(m) \le 0$ and $G_n'(m)<\infty$ for all $n\ge 0$. Back to the moments' language, we obtain $\E(X_n \, m^{X_n})<\infty$ and $\E (m^{X_n}) \ge (m-1) \E(X_n \, m^{X_n})$. The latter expression, by the FKG inequality, is greater than or equal to  $(m-1) \E(X_n)\, \E(m^{X_n})$. Therefore, $\E(X_n) \le \frac{1}{m-1}$, for all $n\ge 0$. By definition, we get $F_\infty=0$.\qed

\bigskip

We now turn to the lower bound for the free energy.

\medskip

\begin{lemma}
\label{l:lb1}

If $\E (m^{X_0})<(m-1)\E(X_0 \, m^{X_0})<\infty$, then there exists $s\in (1, \, m)$
such that
 $$
    (s-1) \E(X_n \, s^{X_n}) - \E (s^{X_n}) \to \infty,
    \qquad
    n\to \infty.
 $$

\end{lemma}

\medskip

\noindent {\it Proof.} Taking $s\in (1, \, m)$ sufficiently close to $m$ we may assure that
$$
   (s-1) \E(X_0 \, s^{X_0}) - \E (s^{X_0}) >0.
$$

\noindent In the language of generating functions, this means
$$
   (s-1)s G_{0}'(s)-G_{0}(s) > 0.
$$

\noindent By \eqref{G10},
\begin{eqnarray*}
    (s-1)s G_{n+1}'(s)-G_{n+1}(s)
 &=& \frac{m}{s} \Big[ s(s-1) G_n'(s) - \frac{s}{m} G_n(s)\Big] G_n(s)^{m-1}
     \\
 &\ge& \frac{m}{s} [ s(s-1) G_n'(s) - G_n(s) ] G_n(s)^{m-1}
     \\
 &\ge& \frac{m}{s} [ s(s-1) G_n'(s) - G_n(s) ],
\end{eqnarray*}

\noindent where at the last step we used $G_n(s) \ge 1$ for all $s>1$.
By induction, we obtain
$$
    (s-1)s G_{n}'(s)-G_{n}(s) \ge   \Big(\frac{m}{s}\Big)^{\! n}  [(s-1)s G_{0}'(s)-G_{0}(s)] \to \infty
$$

\noindent which is precisely equivalent to the claim of our lemma.\qed

\medskip

\begin{lemma}
\label{l:lb2}

 If $F_\infty =0$, then $\sup_{n\ge 0} \E (X_n \, s^{X_n}) < \infty$ for all $s\in (0, \, m)$.

\end{lemma}

\medskip

\noindent {\it Proof.} Let $k\ge 1$ and $n\ge 0$. Clearly,
$$
   X_{n+k} \ge \sum_{i=1}^{m^k} {\bf 1}_{\{ X_n^{(i)} \ge k+1 \} },
$$

\noindent where, as before, $X_n^{(i)}$, $i\ge 1$, are independent copies of $X_n$. It follows that
$$
   \E(X_{n+k})
   \ge
   m^k \, \P(X_n \ge k+1) \, .
$$

\noindent Suppose $F_\infty=0$. Then by \eqref{F}, $\E(X_{n+k}) \le \frac{1}{m-1}$ for all $n\ge 0$ and $k\ge 1$; hence
$$
   \P(X_n \ge k+1) \le \frac{1}{(m-1)m^k}
$$

\noindent for all $n\ge 0$ and $k\ge 1$, implying the assertion of our lemma.\qed

\bigskip

\noindent {\it Proof of Theorem \ref{t:collet-et-al}: Lower bound.} Assume first $\E (m^{X_0})<(m-1)\E(X_0 \, m^{X_0})<\infty$. Let us prove that $F_\infty > 0$.

By Lemma \ref{l:lb1}, there exists $s\in (1, \, m)$ such that
\begin{equation}
    \label{pf1}
    \E(X_n \, s^{X_n}) \to \infty,
    \qquad
    n\to \infty.
\end{equation}

\noindent If $F_\infty=0$, then by Lemma \ref{l:lb2} we would have $\sup_{n\ge 0} \E(X_n \, s^{X_n}) < \infty$, contradicting \eqref{pf1}.

Consider now the remaining case $\E(X_0 \, m^{X_0}) =\infty$. For any $k\ge 1$ let $\widetilde{X}_{0,k} := \min\{ X_0, \, k\}$ be the trimmed version of $X_0$.
By choosing $k$ sufficiently large, one obtains
$$
  (m-1) \E(\widetilde{X}_{0,k} \, m^{\widetilde{X}_{0,k}}) > \E(m^{\widetilde{X}_{0,k}}).
$$

\noindent The just proved part of our theorem asserts that the free energy associated with $\widetilde{X}_{0,k}$ is positive, and a fortiori $F_\infty>0$ in this case.\qed

\section{Proof of Theorem \ref{t:G(Xn)->1}}
\label{s:pf_Theorem_E(Xn)->0}

We prove Theorem \ref{t:G(Xn)->1} in this section by means of Proposition \ref{p:Pi<}. Write as before $G_n(s):= \E (s^{X_n})$.

\medskip

\begin{lemma}
\label{l:elementaire}

Assume $\eqref{variete-critique}$.

{\rm (i)} For any $n\ge 0$, $s\mapsto \frac{G_n(s)^{m-1}}{s}$ is non-increasing
on $[1, \, m]$. In particular, we have $\frac{G_n(s)^{m-1}}{s} \ge \frac{G_n(m)^{m-1}}{m}$ for $n\ge 0$ and $s\in [1, \, m]$.

{\rm (ii)} We have $\sup_{n\ge 0} G_n (m) \le m^{1/(m-1)}$.

{\rm (iii)} We have $\inf_{n\ge 0} \P(X_n=0)\ge \frac{m-2}{m-1}$.

\end{lemma}

\medskip

\noindent {\it Proof.} (i) Since $s\mapsto \frac{s\, G_n'(s)}{G_n(s)}$ is non-decreasing on $[1, \, \infty)$ (this is a general property of moment generating functions, and has nothing to do with assumption \eqref{variete-critique}), we have, for $s\in [1, \, m]$,
$$
     (m-1)\frac{s\, G_n'(s)}{G_n(s)} -1 \le (m-1)\frac{m\, G_n'(m)}{G_n(m)} -1 =0
$$

\noindent by \eqref{variete-critique-Gn}. This implies that
$$
    \frac{\!\d}{\! \d s} (\frac{G_n(s)^{m-1}}{s})
    =  [ (m-1)\frac{s\, G_n'(s)}{G_n(s)} -1 ] \frac{G_n(s)^{m-1}}{s^2}
    \le 0
$$

\noindent for $s\in [1, \, m]$; hence $s\mapsto \frac{G_n(s)^{m-1}}{s}$ is non-increasing on $[1, \, m]$.

(ii) By (i), $\frac{G_n(m)^{m-1}}{m} \le G_n(1)^{m-1} =1$, so $G_n(m)^{m-1} \le m$.

(iii) {F}rom \eqref{variete-critique}, we get
\begin{eqnarray*}
    \P(X_n=0) &=&  \E(m^{X_n}) -\E(m^{X_n}\, {\bf 1}_{\{ X_n\ge 1\} })
    \\
              &=& (m-1) \E( X_n \, m^{X_n}\, {\bf 1}_{\{ X_n\ge 1\} })
                  -\E(m^{X_n}\, {\bf 1}_{\{ X_n\ge 1\} })
    \\
              &\ge& (m-2) \E(m^{X_n}\, {\bf 1}_{\{ X_n\ge 1\} }).
\end{eqnarray*}

\noindent This implies $(m-1)\P(X_n=0)\ge(m-2)\E(m^{X_n})\ge m-2$, as desired.\qed

\bigskip

\noindent {\it Proof of Theorem \ref{t:G(Xn)->1}.} Assume \eqref{variete-critique}. Let, for $n\ge 0$,
$$
   \varepsilon_n
   :=
   G_n(m) -1
   >
   0\, .
$$

\noindent The proof of Theorem \ref{t:G(Xn)->1} consists of showing that $\varepsilon_n \to 0$.

By \eqref{G0}, $G_{n+1}(s) \le \frac1s \, G_n(s)^m + (1-\frac1s)$. In particular, taking $s=m$ gives that
\begin{equation}
    \varepsilon_{n+1}
    \le
    \frac{(1+\varepsilon_n)^m -1}{m} ,
    \qquad
    n\ge 0\, .
    \label{epsilon(n+1)<epsilon(n)}
\end{equation}

We will now use the following elementary inequality
\begin{equation}
    \frac{m}{\log (1+my)}
    <
    \frac{1}{\log (1+y)} + \frac{m-1}{2} \, ,
    \qquad
    y>0\, .
    \label{ineg_log}
\end{equation}

\noindent Indeed, the function $h(x):=\tfrac{x}{\log(1+x)}$ satisfies $h'(x)\le \tfrac{1}{2}$ for all $x>0$. Therefore,
$$
   \frac{m}{\log(1+my)} -  \frac{1}{\log(1+y)}
   = y^{-1} (h(my)-h(y)) \le  \frac{m-1}{2} \, .
$$

Let now $n>k\ge 0$. Since $\varepsilon_{n-k} \le \frac{(1+\varepsilon_{n-k-1})^m -1}{m}$ (see \eqref{epsilon(n+1)<epsilon(n)}), we have
$$
   \frac{1}{\log (1+\varepsilon_{n-k-1})}
   \le
   \frac{m}{\log (1+m \varepsilon_{n-k})}
   <
   \frac{1}{\log (1+\varepsilon_{n-k})} + \frac{m-1}{2} \, ,
$$

\noindent the last inequality being a consequence of \eqref{ineg_log}. Iterating
the inequality yields that for all integers $0\le k\le n$,
$$
     \frac{1}{\log (1+\varepsilon_{n-k})}
     \le
     \frac{1}{\log (1+\varepsilon_n)} + \frac{m-1}{2}\, k
     \le
     \frac{1}{\varepsilon_n} + 1 + \frac{m-1}{2}\, k \, ,
$$

\noindent the second inequality following from the inequality $\frac{1}{\log (1+x)} < \frac1x + 1$ (for $x>0$). As a consequence, for integers $n> j\ge 0$,
$$
     \sum_{k=0}^{n-j-1} \log (1+\varepsilon_{n-k})
     \ge
     \sum_{k=0}^{n-j-1} \frac{1}{\frac{1}{\varepsilon_n} + 1 + \frac{m-1}{2}\, k}
     \ge
     \int_0^{n-j} \frac{\d x}{\frac{1}{\varepsilon_n} + 1 + \frac{m-1}{2}\, x} \, .
$$

\noindent The sum on the left-hand side is $\sum_{i=j+1}^n \log (1+\varepsilon_i)$, whereas the integral on the right-hand side is equal to
$$
     \frac{2}{m-1} \log \Big( \frac{\frac{1}{\varepsilon_n} + 1 + \frac{m-1}{2}\,
     (n-j)}{\frac{1}{\varepsilon_n} + 1} \Big)
     \ge
     \frac{2}{m-1} \log ( c_5 \, (n-j) \varepsilon_n),
$$

\noindent for some constant $c_5>0$ whose value depends only on $m$ (recalling from Lemma \ref{l:elementaire}~(ii) that $\varepsilon_n \le m^{1/(m-1)}-1$). This yields that for $n>j\ge 0$,
\begin{equation}
    \prod_{i=j+1}^n (1+\varepsilon_i)^{(m-1)/2}
    \ge
    c_5 \, (n-j) \varepsilon_n\, .
    \label{epsilon_lb}
\end{equation}

\noindent Replacing the pair $(n, \, j)$ by $(j, \, 0)$, we also have
$$
    \prod_{i=1}^j (1+\varepsilon_i)^{(m-1)/2} \ge c_5 \,j \varepsilon_j
$$

\noindent for $j\ge 1$. Therefore,
$$
    c_5^2 \, j(n-j) \, \varepsilon_j\, \varepsilon_n
    \le
    \prod_{i=1}^n (1+\varepsilon_i)^{(m-1)/2},
$$

\noindent which is bounded by $\frac{c_3^{1/2}}{[\P(X_0=0)]^{1/2}} \, n$ (see Proposition \ref{p:Pi<}). Consequently, with $c_6 := \frac{c_3^{1/2}}{c_5^2}$, we have $\varepsilon_j \, \varepsilon_n \le \frac{c_6}{[\P(X_0=0)]^{1/2}} \, \frac{n}{j (n-j)}$ for $n > j \ge 1$. In particular,
\begin{equation}
    \label{sup_epsilon}
    \varepsilon_j \, \sup_{n\ge 2j} \varepsilon_n
    \le
    \frac{2c_6}{[\P(X_0=0)]^{1/2}} \, \frac1j\, ,
    \qquad
    j \ge 1\, .
\end{equation}

\noindent This yields
$$
   (\limsup_{j\to\infty} \varepsilon_j)^2
   =  \limsup_{j\to\infty} \varepsilon_j\ \cdot \lim_{j\to\infty}  \sup_{n\ge 2j} \varepsilon_n
   =  \limsup_{j\to\infty} \Big( \varepsilon_j \, \sup_{n\ge 2j} \varepsilon_n \Big)
   =0,
$$

\noindent i.e., $\varepsilon_j \to 0$.\qed

\section{Around Conjectures \ref{conj:P(Xn=0)} and \ref{conj:G(Xn)}b}
\label{s:conj_P(Xn=0)}

Let
$$
   \varepsilon_n :=  G_n(m)-1.
$$

By Theorem \ref{t:G(Xn)->1}, $\varepsilon_n \to 0$, $n\to \infty$. Propositions \ref{p:Pi<} and \ref{p:Pi>} together say that $\prod_{i=0}^n (1+\varepsilon_i)$ is of order of magnitude $n^{2/(m-1)}$ when $n$ is large. So if $n\mapsto \varepsilon_n$ were sufficiently regular, we would have
\begin{equation}  \label{G(X_n)(m)-1:reapparition}
    \varepsilon_n
    \, \sim \,
    \frac{2}{m-1} \, \frac1n \, ,
    \qquad
    n\to \infty \, ,
\end{equation}

\noindent This is \eqref{G(X_n)(m)-1} in Conjecture \ref{conj:G(Xn)}b.

{F}rom the relation \eqref{G0} we obtain, with $s=m$,
$$
    1+ \varepsilon_{n+1}
    =
    \frac1m \, (1+\varepsilon_n)^m
    +
    (1-\frac1m) \, G_n(0)^m .
$$

\noindent Since $G_n(0) = 1- \P (X_n \not= 0)$, whereas $\varepsilon_n\to 0$, this implies that
\begin{equation}
    \P (X_n \not= 0)
    \sim
    \frac{\varepsilon_n - \varepsilon_{n+1}}{m-1}
    +
    \frac12 \, \varepsilon_n^2 ,
    \qquad
    n\to \infty.
    \label{proba_decollage2}
\end{equation}

Let us look back at \eqref{G(X_n)(m)-1:reapparition}. If we were able to show that $n^2 (\varepsilon_n - \frac{2}{m-1} \, \frac1n)$ admits a finite limit when $n\to \infty$, \eqref{proba_decollage2} would give an affirmative answer to Conjecture \ref{conj:P(Xn=0)}.

\section{About Conjectures \ref{conj:cvg_faible}, \ref{conj:E(Xn)} and \ref{conj:G(Xn)}a}
\label{s:conj_E(Xn)}

Let us first look at
Conjecture \ref{conj:G(Xn)}a. Theorem \ref{t:G(Xn)->1} yields $\P(X_n \not= 0) \to 0$. Therefore,
$$
    G_n(0)^m = [1- \P(X_n \not= 0)]^m = 1- m \, \P(X_n \not= 0) + o(\P(X_n \not= 0)),
    \quad
    n\to \infty.
$$

\noindent Using this fact and the identity \eqref{G0}, we obtain
\begin{eqnarray*}
    \frac{G_{n+1} (s)-1}{\P(X_n \not= 0)}
 &=&
    \frac{s^{-1}(G_n(s)^m-1) +(1-s^{-1}) (G_n(0)^m-1)   }
         {\P(X_n \not= 0)}
    \\
 &=&
    \frac{s^{-1}(G_n(s)^m -1)}{\P(X_n \not= 0)}
    -
    m(1-s^{-1})
    +
    o(1),
    \qquad
    n\to \infty\, .
\end{eqnarray*}

Suppose we were able to prove that $\P(X_{n+1} \not= 0) \sim \P(X_n \not= 0)$
(which would be a consequence of Conjecture \ref{conj:P(Xn=0)}), and that
\begin{equation}
    \frac{G_n (s)-1}{\P(X_n \not= 0)} \to H(s),
    \label{cvg_etroite}
\end{equation}

\noindent for $s\in (0,m)$, with some real-valued measurable function $H(\, \cdot \, )$. Then
$$
   H(s)
   =
   s^{-1} \, m\, H(s)
   -
   m(1-s^{-1}) ,
$$

\noindent which would lead to
$$
   H(s)
   =
   \frac{m(s-1)}{m-s},
$$

\noindent for $s\in (0, \, m)$. This, in view of Conjecture \ref{conj:P(Xn=0)}, is what have led us to 
Conjecture \ref{conj:G(Xn)}a.

To see why Conjecture \ref{conj:cvg_faible} would be true, let us note from \eqref{cvg_etroite} and $H(s)=\frac{m(s-1)}{m-s}$ that
\begin{equation}
    \E(s^{X_n} \, | \, X_n \not= 0)
    \to
    \frac{(m-1)s}{m-s},
    \qquad
    s\in (0, \, m)\, .
    \label{cvg_etroite_2}
\end{equation}

\noindent On the right-hand side, $\frac{(m-1)s}{m-s} = \E(s^{Y_\infty})$ if
$Y_\infty$ is such that $\P(Y_\infty = k) = \frac{m-1}{m^k}$ for integers
$k\ge 1$. In words, $Y_\infty$ has a geometric distribution with parameter
$\frac1m$. This would give a proof of Conjecture \ref{conj:cvg_faible}.

Conjecture \ref{conj:E(Xn)} would be an immediate consequence of Conjecture \ref{conj:P(Xn=0)} and \eqref{cvg_etroite_2}.

We end this section with an elementary argument showing that our prediction for $\E(X_n)$ in Conjecture \ref{conj:E(Xn)} would be a consequence of Conjecture \ref{conj:P(Xn=0)} (without using the conjectured convergence in \eqref{cvg_etroite_2}). 

\begin{lemma}
\label{l:proba_decollage}

 Assume $\E(X_0)<\infty$. For all $n\ge 0$, we have $\E(X_n)<\infty$, and
 $$
 \P (X_n = 0)
 =
 [1 - m\, \E(X_n) + \E(X_{n+1})]^{1/m} .
 $$

\end{lemma}

\medskip

\noindent {\it Proof.} Taking $s=1$ in \eqref{G1} gives that for $n\ge 0$, $\E(X_n)<\infty$ and
$$
    \E(X_{n+1})
    =
    m\, \E(X_n) -1 + G_n(0)^m \, ,
$$

\noindent which implies the lemma by noting that $G_n(0) = \P (X_n = 0)$.\qed

\bigskip

If Conjecture \ref{conj:P(Xn=0)} were true, then it would follow from Lemma \ref{l:proba_decollage} (and Theorem \ref{t:G(Xn)->1} which guarantees $\E(X_n)\to 0$) that
$$
\E(X_n) - \frac{\E(X_{n+1})}{m}
\, \sim \,
\frac{4}{(m-1)^2} \, \frac{1}{n^2} \, ,
\qquad
n\to \infty\, .
$$

\noindent Applying Lemma \ref{l:suite_decalee} below to $a_n := \E (X_n)$,
$\lambda := \frac1m$ and $b:=2$, this would yield the prediction for $\E(X_n)$ in Conjecture \ref{conj:E(Xn)}.

\begin{lemma}
\label{l:suite_decalee}

 Let $0<\lambda<1$ and $b>0$. Let $(a_n, \, n\ge 1)$ be a bounded sequence such that
 $$
 \lim_{n\to \infty} \, n^b(a_n - \lambda a_{n+1}) = x \, ,
 $$

 \noindent for some real $x$. Then
 $$
    \lim_{n\to \infty} \, n^b a_n
    =
    \frac{x}{1-\lambda} \, .
 $$

\end{lemma}

\medskip

\noindent {\it Proof.} Let $p_n:= a_n- \lambda a_{n+1}$. Then $ \lim_{n\to\infty}  n^b p_n=x$ and
$$
   M:= \sup_{j\ge 1} \, ( j^b |p_j| ) <\infty.
$$

\noindent In particular, the sequence $(p_n)$ is bounded.

By iterating the identity $a_n=p_n + \lambda a_{n+1}$, we get a series representation
$$
   a_n = \sum_{k=0}^\infty p_{n+k} \lambda^k,
$$

\noindent where the series converges because $(p_n)$ is bounded. We may apply the dominated convergence theorem to the identity
$$
  n^b a_n = \sum_{k=0}^\infty n^b p_{n+k} \lambda^k,
$$

\noindent because for every $k\ge 0$ we have $n^b p_{n+k}\to x$ by assumption, and the dominating summable majorant $(M_k, k\ge 0)$ is given by $M_k:= M \lambda^k$. We arrive at
$$
   \lim_{n\to\infty}  n^b a_n = x \sum_{k=0}^\infty \lambda^k =\frac{x}{1-\lambda}\, ,
$$

\noindent as required. \qed

\section{Proofs of Propositions \ref{p:Pi<} and \ref{p:Pi>}}
\label{s:prop_proofs}

\subsection{Proof of Proposition \ref{p:Pi<}}

Let us define, for $n\ge 0$ and $s\in [0, \, m)$,
$$
   \Delta_n(s)
   :=
   [G_n(s) - s(s-1)G_n'(s)]
   -
   \frac{(m-1)(m-s)}{m} \, [2sG_n'(s) + s^2 G_n''(s)]\, .
$$

\noindent Then by \eqref{G10}, \eqref{G1} and \eqref{G2},
\begin{equation}
    \label{Delta(n+1)}
    \Delta_{n+1}(s)
    =
    \frac{m}{s}\, \Delta_n(s) \, G_n(s)^{m-1}
    -
    \frac{m-s}{s}\, [(m-1)sG_n'(s)-G_n(s)]^2  G_n(s)^{m-2} .
\end{equation}

\medskip

\begin{lemma}
 \label{l:Dg}

 Assume $\eqref{variete-critique}$. Then $\Delta_n(s)\in [0,\, 1]$ for $n\ge 0$ and $s\in [0, \, m)$.

\end{lemma}

\medskip

\noindent {\it Proof.} By the definition of $\Delta_n$ , for $s\in [0, \, m)$,
$$
   \Delta_n'(s)
   =
   -\frac{m-s}{m} \, [ 2(m-2)G_n'(s) + (4m-5) sG_n''(s) + (m-1) s^2G_n'''(s) ] ,
$$

\noindent which is non-positive. Hence $\Delta_n$ is non-increasing on $(0, \, m)$. Since $\Delta_n(0) = G_n(0) = \P(X_n=0) \le 1$, whereas under assumption \eqref{variete-critique}, it is easily checked that $\lim_{s\to m-} \Delta_n(s) = 0$, the lemma follows.\qed

\medskip

\begin{lemma}
 \label{l:squareD}

 Assume $\eqref{variete-critique}$. For all $n\ge 0$,
 $$
 [ (m-1)s G_n'(s)-G_n(s) ]^2
 \le
2 G_n(0)\Delta_n(s),
 \qquad
 s\in [0, \, m).
 $$
\end{lemma}

\medskip

\noindent {\it Proof.} By using \eqref{variete-critique-Gn} and writing $x := \frac{s}{m}\in [0, \, 1)$ for brevity, we have
$$
   G_n(s)
   =
   (m-1)m G_n'(m) - G_n(m) + G_n(s)
   =
   \sum_{k=1}^\infty m^k(km-k-1+x^k) \, \P (X_n=k) \, .
$$

\noindent [In particular, $G_n(0) = \sum_{k=1}^\infty m^k(km-k-1) \, \P (X_n=k)$.] Furthermore,
\begin{eqnarray*}
    G_n(s) - (m-1)s G_n'(s)
 &=& \sum_{k=1}^\infty  m^k [km-k-1+x^k-(m-1)kx^k]\, \P (X_n=k)
    \\
 &=& \sum_{k=1}^\infty m^k(km-k-1)(1-x^k)\, \P (X_n=k) \, .
\end{eqnarray*}

\noindent Consequently,
\begin{equation}
    \Delta_n(s)
    =
    \sum_{k=1}^\infty m^k(km-k-1)(1-(k+1)x^k+kx^{k+1}) \, \P (X_n=k) \, .
    \label{e:2157}
\end{equation}

\noindent By the Cauchy--Schwarz inequality, $(\sum_{k=1}^\infty a_k q_k)^2 \le (\sum_{k=1}^\infty q_k)(\sum_{k=1}^\infty a_k^2 q_k)$, with $a_k:=1-x^k$ and $q_k:= m^k(km-k-1) \P(X_n=k)$, the proof of the lemma is reduced to showing the following: for $x\in [0,1]$ and $k\ge 1$,
\begin{equation}
    \label{ineg_xkm}
    (1-x^k)^2 \le 2(1-(k+1)x^k+kx^{k+1}).
\end{equation}

\noindent This can be rewritten as $2 k\, x^k(1-x)\le 1-x^{2k}$, which is proved by
$$
  \frac{1-x^{2k}}{1-x}= 1+x+...+x^{2k-1} \ge 2k x^k
$$

\noindent (using that $2x^k\le  2x^{k-1/2} \le x^{k+\ell}+ x^{k-\ell-1}$ for $0\le \ell < k$).\qed

\medskip

\noindent {\it Proof of Proposition \ref{p:Pi<}.}
By \eqref{Delta(n+1)} and Lemma \ref{l:squareD}, for $s\in [0, \, m)$,
\begin{eqnarray*}
    \Delta_{n+1}(s)
 &\ge&\frac{m}{s} \,  \Delta_n(s)\, G_n(s)^{m-1}
    -
    \frac{m-s}{s}  \, 2 G_n(0) \Delta_n(s) \, G_n(s)^{m-2}
    \\
 &\ge&\frac{m}{s} \, \Delta_n(s) \, G_n(s)^{m-1}
    -
    \frac{2(m-s)}{s} \, \Delta_n(s) \, G_n(s)^{m-1}
    \\
 &=& \frac{2s-m}{s} \, \Delta_n(s) \, G_n(s)^{m-1} .
\end{eqnarray*}

\noindent Hence for any $s\in (\frac{m}{2},\, m)$,
$$
\Delta_{n+1}(s)
\ge
\Delta_0(s) (2s-m)^{n+1} \prod_{i=0}^n \frac{G_i(s)^{m-1}}{s}
\ge
\Delta_0(s) (2s-m )^{n+1} \prod_{i=0}^n \frac{G_i(m)^{m-1}}{m} ,
$$

\noindent where we used Lemma \ref{l:elementaire}~(i) in the last inequality.
Therefore, for $s\in (\frac{m}{2}, \, m)$,
$$
\prod_{i=0}^n G_i(m)^{m-1}
\le
\Big(\frac{m}{2s-m}\Big)^{n+1} \, \frac{\Delta_{n+1}(s)}{\Delta_0(s)}
\le
\Big(\frac{m}{2s-m}\Big)^{n+1} \frac{1}{\Delta_0(s)},
$$

\noindent the second inequality being a consequence of Lemma \ref{l:Dg}.
Taking $s=m-\frac{m}{n+2}$ gives that with $c_7:= \sup_{n\ge 1} (1+\frac2n)^{n+1}\in (1, \, \infty)$ and all $n\ge 1$,
$$
\prod_{i=0}^n G_i(m)^{m-1}
\le
\frac{c_7}{\Delta_0(m-\frac{m}{n+2})} \, .
$$

\noindent By \eqref{e:2157} and \eqref{ineg_xkm}, with $x:= 1- \frac{1}{n+2}$,
\begin{eqnarray*}
    \Delta_0(m-\frac{m}{n+2})
 &\ge& \frac12
    \sum_{k=1}^\infty m^k(km-k-1)(1-x^k)^2 \, \P (X_0=k)
    \\
 &\ge& \frac12 \, (1-x)^2 \sum_{k=1}^\infty m^k(km-k-1) \, \P (X_0=k) \, .
\end{eqnarray*}

\noindent On the right-hand side, $(1-x)^2 = \frac{1}{(n+2)^2}$, whereas 
\begin{eqnarray*}
    \sum_{k=1}^\infty m^k(km-k-1) \, \P (X_0=k)
 &=& \E \{ [(m-1)X_0-1]m^{X_0} \, {\bf 1}_{\{ X_0 \ge 1\} }\} 
    \\
 &=& \E \{ [(m-1)X_0-1]m^{X_0}\} + \P(X_0=0),
\end{eqnarray*}

\noindent which is $\P(X_0=0)$ because under \eqref{variete-critique}, we have $\E \{ [(m-1)X_0-1]m^{X_0}\}=0$. Consequently,
$$
\Delta_0(m-\frac{m}{n+2}) \ge \frac{\P(X_0=0)}{2(n+2)^2} \, ,
$$

\noindent for all $n\ge 1$. This yields the proposition.\qed

\subsection{Proof of Proposition \ref{p:Pi>}}
\label{subs:proof_Pi>}

Now we start to prepare the proof of Proposition \ref{p:Pi>}. In the case $m=2$, it is proved in Appendix III of \cite{collet-eckmann-glaser-martin}.

\medskip

\begin{lemma}
\label{l:Gn_derivee_deuxieme(m)}

Assume $\eqref{variete-critique}$ and let $n\ge 0$. If $G_n'''(m)<\infty$,
then
\begin{equation}
    \label{Gn_derivee_deuxieme(m)}
    G_n''(m)
    \le
    c_8\, \max\{G_n'''(m)^{1/2},1\} \, ,
\end{equation}

\noindent where $c_8=c_8(m) \in (0, \, \infty)$ is a constant whose value depends only on $m$.

\end{lemma}

\medskip

\noindent {\it Proof.} According to Lemma \ref{l:elementaire}~(ii), $G_n(m) \le m^{1/(m-1)}$. By plugging this into \eqref{variete-critique}, we have
$G_n'(m) \le \frac{m^{1/(m-1)}}{m(m-1)} =: c_9$. The function $G_n'(\, \cdot\, )$ being convex, we have
$$
   G_n''(s)
   \le
   \frac{G_n'(m)-G_n'(s)}{m-s}
   \le
   \frac{G_n'(m)}{m-s}
   \le
   \frac{c_9}{m-s},
   \qquad
   s\in [0, \, m) \, .
$$

\noindent By the convexity of $G_n''(\, \cdot\, )$, this implies that, for $s\in [0, \, m)$,
$$
   G_n''(m)
   \le
   G_n''(s)
   +
   (m-s) \, G_n'''(m)
   \le
   \frac{c_9}{m-s}
   +
   (m-s) \, G_n'''(m) \, .
$$

\noindent The lemma follows by taking $s:= m - \frac{1}{\max\{G_n'''(m)^{1/2},1\}}$.\qed

\medskip

\begin{lemma}
\label{l:G-1>}

 Assume $\eqref{variete-critique}$ and let $n\ge 0$. If $G_n'''(m)<\infty$, then
 \begin{equation}  \label{G-1>}
    G_n(m)-1
    \ge
    \frac{c_{10}}{ \max\{G_n'''(m)^{1/2},1\} } \, ,
 \end{equation}

 \noindent where $c_{10}=c_{10}(m) \in (0, \, \infty)$ is a constant that
 does not depend on $n$.

\end{lemma}

\medskip

\noindent {\it Proof.} If $m=2$ and $\P(X_n=0)\le \tfrac12$, then $G_n(2)-1\ge \P(X_n\ge 1)\ge \tfrac12$.  If $m=2$ and $\P(X_n=0)> \tfrac12$, the equation (AIII.5) of \cite{collet-eckmann-glaser-martin} says that $G_n''(2) (G_n(2)-1) > \tfrac14 \, [\P (X_n=0)]^2>\tfrac1{16}$. So the lemma (in case $m=2$) follows from Lemma \ref{l:Gn_derivee_deuxieme(m)}.

In the rest of the proof, we assume $m\ge 3$.

Write $\eqref{variete-critique}$ as
$$
    \sum_{k=0}^\infty ((m-1)k-1) m^k \, \P(X_n=k) = 0.
$$

\noindent It follows that
\begin{eqnarray*}
    \P(X_n=0) - m(m-2)\, \P(X_n=1)
 &=& \sum_{k=2}^\infty ((m-1)k-1) m^k \, \P(X_n=k)
    \\
 &\le& (m-1)\sum_{k=2}^\infty k m^k \, \P(X_n=k) .
\end{eqnarray*}

\noindent Writing
\begin{eqnarray*}
    G_n''(m)
 &=& \frac{1}{m^2} \, \sum_{k=2}^\infty k(k-1)m^k \, \P(X_n=k)
    \ge
    \frac{1}{2m^2} \, \sum_{k=2}^\infty k^2\, m^k \, \P(X_n=k),
    \\
 G_n(m) -1
 &=& \sum_{k=1}^\infty (m^k-1)\, \P(X_n=k)
    \ge
    (1- \frac1m) \sum_{k=1}^\infty m^k \, \P(X_n=k) \, ,
\end{eqnarray*}

\noindent it follows from the Cauchy--Schwarz inequality that
$$
\P(X_n=0)
-
m(m-2) \, \P(X_n=1)
\le
\Big( 2m^3 (m-1) \, G_n''(m) [G_n(m) -1]\Big)^{1/2} \, .
$$

\noindent Write $\varepsilon_n := G_n(m)-1$ as before. We have $G_n(m)\ge 1+(m-1)\P(X_n=1)$, i.e.,
$$
   m(m-2) \, \P(X_n=1) \le \frac{m(m-2)\, \varepsilon_n}{m-1}.
$$

\noindent On the other hand, $\P(X_n=0) \ge \frac{m-2}{m-1}$ by Lemma \ref{l:elementaire}~(iii). Therefore,
$$
   \frac{m-2}{m-1} - \frac{m(m-2)\, \varepsilon_n}{m-1}
   \le
   (2m^3 (m-1))^{1/2}\, G_n''(m)^{1/2} \, \varepsilon_n^{1/2} ,
$$

\noindent which yields
\begin{equation}
    G_n(m) -1
    =
    \varepsilon_n
    \ge
    \frac{c_{11}}{\max\{G_n''(m),1\}},
    \label{e:loweps}
\end{equation}

\noindent for some constant $c_{11} = c_{11}(m) \in (0, \, \infty)$ and all $n$. The lemma (in case $m\ge 3$) follows from Lemma \ref{l:Gn_derivee_deuxieme(m)}.\qed

\medskip

\noindent {\it Proof of Proposition \ref{p:Pi>}.} Write $b_j= \prod_{i=0}^{j}G_i(m)^{m-1}$. By \eqref{D2(m):lemme},
$$
   \DD_n(m)=\DD_0(m) b_{n-1} \, .
$$

\noindent By definition \eqref{D2(m)}, we have $\DD_n(m) \ge m(m-1)\, G_n'''(m)$ and $\frac{\DD_0(m)}{m(m-1)} \le c_{12} \Lambda_0$ for some $c_{12} = c_{12}(m) \in (0, \, \infty)$, where $\Lambda_0 := \E (X_0^3 m^{X_0})$. Therefore,
\begin{equation}
    \label{e:G'''low}
    G_n'''(m)
    \le
    c_{12} \, \Lambda_0\, b_{n-1} \, .
\end{equation}

\noindent We know that $b_{n-1}\ge 1$, and that $\Lambda_0 \ge \E (X_0 m^{X_0}) = \frac{1}{m-1} \E (m^{X_0}) \ge \frac{1}{m-1}$. So, with $c_{13} := \max\{ c_{12}, \, m-1\}$, we have
$$
   \max\{ G_n'''(m), \, 1\}
   \le
   c_{13} \, \Lambda_0\, b_{n-1} \, ,
$$

\noindent On the other hand, by Lemma \ref{l:G-1>}, $\max\{ G_n'''(m), \, 1\} \ge (\frac{c_{10}}{G_n(m)-1})^2$. Therefore, with $c_{14} := \frac{c_{10}}{c_{13}^{1/2}}$, we have $G_n(m) \ge 1+ \frac{c_{14} }{\Lambda_0^{1/2} b_{n-1}^{1/2}}$, i.e.,
$$
   \frac{b_n^{1/(m-1)}}{b_{n-1}^{1/(m-1)}}
   \ge
   1+ \frac{c_{14}}{\Lambda_0^{1/2} b_{n-1}^{1/2}}\, .
$$

\noindent Since $\frac{c_{14}}{\Lambda_0^{1/2} b_{n-1}^{1/2}} \le \frac{c_{14}}{\Lambda_0^{1/2}} \le (m-1)^{1/2}c_{14}$, this yields that for some $c_{15}= c_{15}(m) \in (0, \, \infty)$ depending only on $m$,
$$
    \frac{b_n^{1/2}}{b_{n-1}^{1/2}}
    \ge
    1+ \frac{c_{15} }{\Lambda_0^{1/2} b_{n-1}^{1/2}}\, .
$$

\noindent Hence $b_n^{1/2} - b_{n-1}^{1/2} \ge \frac{c_{15}}{\Lambda_0^{1/2}} $. Thus $b_n^{1/2}\ge \frac{c_{15}}{\Lambda_0^{1/2}} n$, as desired.
\qed

\section{Proof of Theorem \ref{t:G(Xn)(m)-1}}
\label{s:pf_Theorem_G(Xn)(m)-1}

Assume $\eqref{variete-critique}$. Write as before $G_n(m) := \E(m^{X_n})$
and $\varepsilon_n := G_n(m) -1$. Recall from \eqref{epsilon_lb} that for
integers $n>j\ge 0$,
$$
    \prod_{i=j+1}^n (1+\varepsilon_i)^{(m-1)/2} \ge c_5 \, (n-j) \varepsilon_n.
$$

\noindent Taking $j:= \lfloor \frac{n}{2} \rfloor$ gives that for $n\ge 1$,
$$
\prod_{i=\lfloor \frac{n}{2} \rfloor+1}^n (1+\varepsilon_i)^{(m-1)/2}
\ge
\frac{c_5}{2} \, n \varepsilon_n\, .
$$

\noindent By Proposition \ref{p:Pi<},
$$
   \prod_{i=0}^n (1+\varepsilon_i)^{(m-1)/2} \le \frac{c_3^{1/2}}{\P(X_0=0)^{1/2}} \, n,
$$

\noindent whereas by Proposition \ref{p:Pi>}, under the assumption $\Lambda_0:= \E(X_0^3 \, m^{X_0})<\infty$,
$$
   \prod_{i=0}^{\lfloor \frac{n}{2} \rfloor} (1+\varepsilon_i)^{(m-1)/2}
   \ge
   \frac{c_{16}}{\Lambda_0^{1/2}} \, n
$$

\noindent (for some constant $c_{16} = c_{16}(m)>0$ depending only on $m$; noting that the case $n=1$ is trivial because $\Lambda_0 \ge \frac{1}{m-1}$). Hence
$$
    \prod_{i=\lfloor \frac{n}{2} \rfloor+1}^n (1+\varepsilon_i)^{(m-1)/2}
    \le
    \frac{c_3^{1/2} \, \Lambda_0^{1/2}}{c_{16} \, \P(X_0=0)^{1/2}} .
$$

\noindent Consequently, $\varepsilon_n \le c_{17} \frac{\Lambda_0^{1/2}}{\P(X_0=0)^{1/2}} \, \frac1n$ with $c_{17} := \frac{2c_3^{1/2}}{c_5c_{16}}$.

To prove the lower bound, we note that under the assumption $\Lambda_0 := \E(X_0^3 \, m^{X_0})<\infty$, we have, by  (\ref{G-1>}),  $\varepsilon_n \ge \frac{c_{10}}{\max\{G_n'''(m)^{1/2},1\} }$; since
$$
     G_n'''(m)
     \le
     c_{12} \Lambda_0 \prod_{i=0}^{n-1} (1+\varepsilon_i)^{m-1}
$$

\noindent (see \eqref{e:G'''low}), which is bounded by $c_{12} \Lambda_0 \, \frac{c_3}{\P (X_0=0)} \, n^2$ according to Proposition \ref{p:Pi<}, this yields that
$$
    \varepsilon_n
    \ge
    \frac{c_{10}}{\max \Big\{ (c_{12}c_3)^{1/2} \frac{\Lambda_0^{1/2}}{[\P (X_0=0)]^{1/2}}\, n, \, 1\Big\}} ,
$$

\noindent as desired.\qed

\section{Some additional rigorous results}
\label{s:harmonic_limit}

Let $(b_i)$ be a numerical sequence. We define its {\it harmonic limit} as
$$
   h\hbox{-}\!\lim_{i\to\infty} b_i
   :=
   \lim_{n\to \infty} \frac1{\log n} \sum_{i=1}^{n} \frac{b_i}i\, ,
$$

\noindent if the latter exists. The $h\hbox{-}\!\limsup$ and  $h\hbox{-}\!\liminf$ are defined in the same way.

\begin{proposition}
\label{p:h-lim_eps}

Assume $\eqref{variete-critique}$. If $\,\E(X_0^3 \, m^{X_0})<\infty$, then
$$
    h\hbox{-}\!\lim_{i\to \infty}[ (\E( m^{X_i})-1)i] = \frac{2}{m-1}\, .
$$

\end{proposition}

This statement is a (much) weaker version of Conjecture \ref{conj:G(Xn)}b and a complement to Theorem \ref{t:G(Xn)(m)-1}. However, it has the advantage that the limiting constant $\tfrac{2}{m-1}$ appears explicitly in the result.

\bigskip

\noindent {\it Proof.} We denote as usual $\varepsilon_i:= G_i(m)-1=\E( m^{X_i})-1$. Rewrite the assertions of Propositions \ref{p:Pi<} and  \ref{p:Pi>} as
$$
  \frac{1}{m-1} \log \frac{c_4}{\Lambda_0} +  \frac{2}{m-1}\, \log n
  \le
  \sum_{i=1}^{n} \log(1+\varepsilon_i)
  \le
  \frac{1}{m-1} \log \frac{c_3}{\P(X_0=0)} +  \frac{2}{m-1}\, \log n.
$$

\noindent Using $\log(1+x) \le x \le \log(1+x)+ x^2/2$ for all $x>0$ we obtain
$$
  \frac{1}{m-1} \log \frac{c_4}{\Lambda_0} +  \frac{2}{m-1}\, \log n
  \le
  \sum_{i=1}^{n} \varepsilon_i
  \le
  \frac{1}{m-1} \log \frac{c_3}{\P(X_0=0)} +  \frac{S}2 +\frac{2}{m-1}\, \log n
$$

\noindent where $S:= \sum_{i=1}^{\infty} \varepsilon_i^2<\infty$ by Theorem \ref{t:G(Xn)(m)-1}. It follows immediately that
\begin{equation}
    \label{sum_epsilon}
    \lim_{n\to \infty} \frac1{\log n} \sum_{i=1}^{n} \frac{\varepsilon_i\, i}i\
    =
    \lim_{n\to \infty} \frac1{\log n} \sum_{i=1}^{n} \varepsilon_i
    =
    \frac{2}{m-1}\, ,
\end{equation}

\noindent proving the proposition.\qed

\medskip

\begin{proposition}
\label{p:h-lim_P}

Assume $\eqref{variete-critique}$. If $\,\E(X_0^3 \, m^{X_0})<\infty$, then
$$
  0 <     h\hbox{-}\!\liminf_{i\to \infty} [ \P(X_i\not=0) i^2]
    \le   h\hbox{-}\!\limsup_{i\to \infty} [\P(X_i\not=0) i^2]
    < \infty\, .
$$
\end{proposition}

This result is a weaker version of Conjecture \ref{conj:P(Xn=0)}, but at least we are able to say something rigorous about the order $i^{-2}$ for $\P(X_i\not=0)$.

\medskip

\noindent {\it Proof.} It is more convenient to use in the calculations the quantity
$$
p_i := m^{-1}(1-G_i(0)^m) \sim \P(X_i\not=0) ,
\qquad
i\to \infty\, .
$$

\noindent Recall that \eqref{G0} with $s=m$ may be written as
$$
     1+ \varepsilon_{i+1} = \frac {1}{m} (1+\varepsilon_i)^m + \Big(1-\frac {1}{m}\Big) (1-m p_i)
$$

\noindent whereas
$$
\varepsilon_{i+1}
=
\varepsilon_i+ \frac{m-1}2\, \varepsilon_i^2 (1+o(1)) - (m-1) p_i,
\qquad
i\to \infty\, ,
$$

\noindent or equivalently, as in \eqref{proba_decollage2},
$$
p_i
=
\frac{\varepsilon_i - \varepsilon_{i+1}}{m-1} + \frac{\varepsilon_i^2}2\, (1+o(1)),
\qquad
i\to \infty\, .
$$

\noindent It follows that, for $n\to \infty$,
$$
  \sum_{i=1}^{n} \frac{p_i\, i^2}i
  =
  \sum_{i=1}^{n} [p_i\, i]
  =
   \frac1{m-1} \Big[  - n  \varepsilon_{n+1}+  \sum_{i=1}^{n} \varepsilon_i \Big]
   + \sum_{i=1}^{n} \frac{\varepsilon_i^2 i }2\, (1+o(1)).
$$

\noindent By using \eqref{sum_epsilon} and Theorem \ref{t:G(Xn)(m)-1}, we immediately obtain
\begin{equation}
    \label{h-linf_pi}
    h\hbox{-}\!\liminf_{i\to \infty} [p_i i^2]
    \ge
    \frac{2}{(m-1)^2}  +   \frac{c_1^2}2 \, \frac{\P(X_0=0)}{\Lambda_0}
    > 0
\end{equation}

\noindent and
\begin{equation}
    \label{h-lsup_pi}
    h\hbox{-}\!\limsup_{i\to \infty} [p_i i^2]
    \le
    \frac{2}{(m-1)^2}  +   \frac{c_2^2}2 \, \frac{\Lambda_0}{\P(X_0=0)}
    <\infty,
\end{equation}

\noindent where $c_1$ and $c_2$ are the constants from Theorem \ref{t:G(Xn)(m)-1}.\qed

\medskip

\begin{proposition}
\label{p:h-lim_E}

 Assume $\eqref{variete-critique}$. If $\,\E(X_0^3 \, m^{X_0})<\infty$, then
 $$
  0 <     h\hbox{-}\!\liminf_{i\to \infty} [ \E(X_i) i^2]
    \le   h\hbox{-}\!\limsup_{i\to \infty} [ \E(X_i) i^2]
    < \infty\, .
 $$

\end{proposition}

This result is a weaker version of Conjecture \ref{conj:E(Xn)}, but gives a rigorous statement about the order $i^{-2}$ for $\E(X_i)$.

\medskip

\noindent {\it Proof.} As in the previous proof, we will use $p_i := m^{-1}(1-G_i(0)^m)$. We already know that
$$
  \E(X_{i+1})= m \E(X_i) -1+ G_n(0)^m =  m \E(X_i) -m p_i,
$$

\noindent whereas
$$
   p_i  =  \E(X_i) - \frac 1m \, \E(X_{i+1}).
$$

\noindent By summing up, we obtain
$$
  \sum_{i=1}^n [p_i\, i]
  =
  \frac{m-1}{m}\sum_{i=2}^n [\E(X_i)\, i]
  +
  \frac{1}{m} \sum_{i=2}^n \E(X_i)
  +
  \E(X_{1}) - \frac  {\E(X_{n+1}) n}{m}.
$$

\noindent Notice that the last term is bounded because by Jensen's inequality,
$$
    \E(X_n) \log m \le \log \E (m^{X_n}) = \log (G_n(m)) \le G_n(m) -1 \le c_2 \, \frac{\Lambda_0^{1/2}}{[\P(X_0=0)]^{1/2}} \, \frac1n ,
$$

\noindent using Theorem \ref{t:G(Xn)(m)-1} at the last step. It follows that
$$
    \frac{m-1}{m}\ \sum_{i=1}^n [\E(X_i)\, i] + O(1)
    \le
    \sum_{i=1}^n [p_i\, i]
    \le
    \sum_{i=1}^n [\E(X_i)\, i] + O(1),
    \qquad
    n\to \infty\, .
$$

\noindent The proof is completed by an application of \eqref{h-linf_pi} and  \eqref{h-lsup_pi}.\qed

\bigskip
\bigskip

\noindent {\bf Acknowledgments.} X.~C. was supported by NSFC grants Nos.~11771286 and 11531001. M.~L.\ was supported by RFBR grant 16-01-00258. Part of the work was carried out when M.~L.\ and Z.~S.\ were visiting, respectively, LPMA Universit\'e Pierre et Marie Curie in June and July 2016, and New York University Shanghai in spring 2016; we are grateful to LPMA and NYUSH for their hospitality.

\baselineskip=14pt

{\scriptsize 

Xinxing Chen, School of Mathematical Sciences, Shanghai Jiaotong University, 200240 Shanghai, China, {\tt chenxinx@sjtu.edu.cn}

\medskip

Bernard Derrida, Coll\`ege de France, 11 place Marcelin Berthelot, F-75231 Paris Cedex 05, France, and Laboratoire de Physique Statistique, \'Ecole Normale Sup\'erieure, Universit\'e Pierre et Marie Curie, Universit\'e Denis Diderot, CNRS, 24 rue Lhomond, F-75231 Paris Cedex 05, France, {\tt derrida@lps.ens.fr}

\medskip

Yueyun Hu, LAGA, Universit\'e Paris XIII, 99 av.~J-B Cl\'ement, F-93430 Villetaneuse, France, {\tt yueyun@math.univ-paris13.fr}

\medskip

Mikhail Lifshits, St.~Petersburg State University, Russia, and MAI, Link\"oping University, Sweden, {\tt mikhail@lifshits.org}

\medskip

Zhan Shi, LPMA, Universit\'e Pierre et Marie Curie, 4 place Jussieu, F-75252 Paris Cedex 05, France, {\tt zhan.shi@upmc.fr}

} 

\end{document}